\documentclass[12pt]{amsart}

\usepackage{amssymb, amsmath, epsfig,graphics}

\def\mapdownright#1{\Big\downarrow\rlap{$\vcenter{\hbox{$\scriptstyle#1$}}$}}

\newtheorem{theorem}{Theorem}[section]

\newtheorem{proposition}[theorem]{Proposition}

\newtheorem{examples}[theorem]{Examples}

\newtheorem{remark}[theorem]{Remark}
\renewcommand{\proof}{\vspace{.05in}
                    \noindent {\sc Proof} \hspace{.05in}}
\newcommand{\ethrm}{\hspace*{\fill}
                      $\Box$
                      \vspace{.1in}}

\newcommand{\ra}{\rightarrow}

\newcommand{\pic}{\mathop{\rm Pic\,}\nolimits}

\renewcommand{\phi}{\varphi}

\newcommand{\ch}{{\rm ch}}

\renewcommand{\Box}{\square}
\renewcommand{\tilde}{\widetilde}

\newcommand{\A}{\mathbb A}

\newcommand{\h}{{\mathfrak h}}

\newcommand{\X}{\mathcal X}

\newcommand{\Y}{\mathcal Y}
\newcommand{\E}{\mathcal E}
\newcommand{\C}{\mathbb C}
\newcommand{\R}{\mathbb R}
\newcommand{\Z}{\mathbb Z}
\renewcommand{\O}{\mathcal O}

\renewcommand{\P}{\mathbb P}
\newcommand{\Q}{\mathbb Q}

\newcommand{\color}[6]{}

\newcommand{\db}[1]{D^b(#1)}
\newcommand{\daut}[1]{{\rm Auteq}(#1)}

\newcommand{\tensor}{\stackrel{{\bf L}}{\otimes}}

\newcommand{\push}{{\bf R}}
\newcommand{\pull}{{\bf L}}
\newcommand{\fmt}{Fourier--Mukai functor} 

\begin{document}
\begin{center}
{\LARGE Enhanced gauge symmetry and braid group actions}

\vspace{0.2in}

{\large Bal\'azs Szendr\H oi}

\vspace{0.15in}

{\large October 2002}

\end{center}

\vspace{0.15in}

{\small\begin{center} {\sc abstract} \end{center}
{\leftskip=30pt \rightskip=30pt
Enhanced gauge symmetry appears in Type II string theory 
(as well as F- and M-theory) compactified on Calabi--Yau manifolds containing
exceptional divisors meeting in Dynkin configurations. It is shown that in many 
such cases, at enhanced symmetry points in moduli a braid group acts on the 
derived category of sheaves of the variety. This braid group covers the Weyl group 
of the enhanced symmetry algebra, which itself acts on the deformation space 
of the variety in a compatible fashion. 
Extensions of this result are given for nontrivial $B$-fields on K3 surfaces, 
explaining physical restrictions on the $B$-field, as well as for elliptic 
fibrations. The present point of view also gives new evidence for the enhanced 
gauge symmetry content in the case of a local $A_{2n}$-configuration in a threefold 
having global $\Z/2$ monodromy. 
\par}}

\vspace{0.15in}

\section*{Introduction} 

The phenomenon that Type II string theory compactified on a Calabi--Yau manifold
can exhibit enhanced gauge symmetry was first observed in the physics literature 
in the context of K3 surfaces~\cite{wit},~\cite{aK3}. The existence of 
non-perturbatively enhanced symmetry algebras is forced by the 
duality between heterotic string theory on~$T^4$ and the Type IIA string on K3, 
since the former obviously has enhanced symmetry at special points in moduli. 
It was found that a K3 surface can have enhanced gauge symmetry if it has
rational double point (ADE) singularities, and the type of the (simply-laced) 
non-abelian Lie algebra that appears precisely matches that of the singularity. 
The argument for non-abelian gauge symmetry was later extended
to Calabi--Yau threefolds in~\cite{aCY} and~\cite{kmp}, for threefolds with
a curve of ADE singularities. In the presence of monodromy, 
non-simply laced Lie algebras can also appear. These symmetries
and the arising representations have also been analyzed in the context 
of M- and F-theory (see~\cite{ims}, \cite{akm} and references therein). 

The purpose of the present paper is to give a mathematical interpretation of 
a ``holomorphic shadow'' of this symmetry. Namely, of the parameters needed to 
specify a string vacuum, I~will only concentrate on the complex structure and
B-field parameters, ignoring the K\"ahler structure. In particular, by moving in
the K\"ahler moduli 
space, I~can resolve the singularities mentioned in the previous 
paragraph, and work with smooth K3 surfaces and Calabi--Yau threefolds, 
containing ADE configurations of rational curves and configurations of ruled
surfaces respectively. The phenomenon that I~will illustrate by several theorems 
is that enhanced gauge symmetry can occur at points in complex moduli when 
the derived category of the corresponding Calabi--Yau manifold has a large set 
of autoequivalences. Moreover, these derived equivalences always satisfy the
relations of a generalized braid group, which covers the Weyl group of the 
enhanced gauge symmetry Lie algebra. When one deforms the complex parameters, these 
autoequivalences deform away to equivalences of derived categories between 
different manifolds; this is always governed by a Weyl group action on the 
deformation space. In particular, one can phrase the results of this paper as saying
that the category of topological D-branes on a Calabi--Yau 
compactification (cf.~\cite{doug})
has an extra braid group worth of symmetries at enhanced gauge symmetry points,
not present at generic points in moduli.

Braid group actions for groups of Type~$A$ (and~$DE$) 
on derived categories were first constructed in~\cite{st}.
In Section~\ref{sec!K3} of the present paper I~will show how to extend
these actions in two dimensions (K3 surfaces) to cover deformations, 
and how this fits into the framework of 
enhanced gauge symmetry. The autoequivalences will be generalized to 
cover deformations with nonzero~$B$-field; in particular I~will derive
the restrictions on the B-field found in~\cite{aK3} by a duality argument. 

Calabi--Yau threefolds, as mentioned before, can exhibit gauge symmetries 
of all~$A\ldots G_2$ types. Corresponding braid group actions are constructed 
in~\cite{sz_br}. I~explain in Section~\ref{sec!CY} the main points
of the construction, refering back to the (easier) surface case. 
I~also give some examples, including an amusing projective example
exhibiting non-trivial monodromy, and make some comments related to the 
interpretation of the actions as enhanced gauge symmetry.
 
The threefolds appearing in this paper represent the simplest case of enhanced
gauge symmetry, that of ``uniform singularities'' or geometrically ruled 
surfaces (no hypermultiplets in physics-speak). In case there are extra rational 
curves in fibers, the mathematics is more complicated (compare for 
example~\cite{ims}); dissident curves can be flopped, 
there are many more autoequivalences and derived equivalences
around, and it appears to be difficult to formulate a clean statement. 
However, for one highly singular situation studied for example 
in~\cite[Section 4]{akm}, the ideas of the present paper are strong enough to 
provide supporting evidence (though alas not a proof) 
for the gauge symmetry content. The argument is spelled out in 
Remark~\ref{specialcase}.  

The paper begins with two introductory sections: Section~\ref{sec!We} 
recalls reflection groups and (generalized) 
braid groups, whereas Section~\ref{sec!fam} deals with
(families of) equivalences of derived categories. The latter section contains a 
statement which may be of independent interest, connecting deformations
of a \fmt\ of a Calabi--Yau variety with its action on cohomology. 
Section~\ref{sec!aff} points out an extension of the results to elliptic 
fibrations and braid groups of affine type which may be interesting from the 
point of view of F-theory, whereas Section~\ref{sec!mir} poses a challenge 
for symplectic geometry via mirror symmetry. 

\vspace{0.1in}

\section{Reflection groups and generalized braid groups} 
\label{sec!We}

A Dynkin diagram~$\Delta$ in this paper means an irreducible finite type diagram 
corresponding to a finite root system~$\Sigma\subset\h_\R$ 
in a real Euclidean inner 
product space~$\h_\R$. It is well known that such diagrams can be of type
$A_n, B_n, C_n, D_n, E_6, E_7, E_8, F_4$ or~$G_2$.
The root system~$\Sigma$ defines a finite reflection group 
$W_\Delta=\langle r_i\rangle$ acting on~$\h_\R$, 
generated by a set of reflections~$r_1,\ldots ,r_n$ indexed by nodes of 
$\Delta$, equivalently by a set of simple roots. As an abstract group, 
\[ W_\Delta\cong\Big\langle r_i\colon i\in{\rm Nodes}(\Delta)\Big\rangle\Big/\Big\langle r_i^2 = 1, (r_i r_j)^{m_{ij}}= 1 \Big\rangle\]
with one relation for every node~$i$ and one for every pair of different nodes~$(i,j)$
with label~$m_{ij}$. The group~$W_\Delta$ also acts 
on the complex vector space~$\h=\h_\R\otimes\C$.

Define the {\it (generalized) braid group} (also called Artin group)~$B_\Delta$ 
by generators and relations as 
\begin{equation}\label{braidrelations} B_\Delta  = \Big\langle R_i\colon i\in{\rm Nodes}(\Delta)\Big\rangle\Big/\Big\langle \underbrace{R_i R_j\ldots}_{m_{ij}} = \underbrace{R_j R_i\ldots}_{m_{ij}}\Big\rangle 
\end{equation}
with one relation for every pair of different nodes~$(i,j)$ of~$\Delta$, the 
{\em braid relation}. 
There is a group homomorphism~$B_\Delta\rightarrow W_\Delta$ sending 
$R_i$ to~$r_i$. As an example, in the familiar case of type~$A_n$ 
the group~$W_\Delta$ is the symmetric group on~$(n+1)$ letters, whereas 
$B_\Delta$ is the classical braid group on~$(n+1)$ strings. 

\section{Families of derived equivalences}
\label{sec!fam}

If~$X$ is a smooth projective variety, let~$\db{X}$ denote 
the bounded derived category of coherent sheaves on~$X$. 
A {\it kernel} (derived correspondence) between smooth projective varieties 
$X_1,X_2$ is an object~$U\in\db{X_1\times X_2}$. Such an object defines  
a functor
\[\Psi^U\colon \db{X_2}\rightarrow\db{X_1}
\]
by 
\[\Psi^U(-)= \push p_{1*}(U\tensor p_2^*(-)), 
\]
with~$p_i\colon X_1\times X_2\ra X_i$ the projections. If~$\Psi^U$ is an equivalence of 
triangulated categories,  then it is called a {\it Fourier--Mukai functor} and
$U$ is said to be {\em invertible}. 

Suppose that~$\pi\colon \X\ra S$ is a smooth family of projective 
varieties over a complex base~$S$. A {\em relative kernel} is a pair~$(U,\phi)$, 
where  
\begin{itemize}
\item~$\phi\colon S\rightarrow S$ is an analytic automorphism, 
giving rise to the fibre product diagram 
\[\begin{array}{ccc}
\X\times_\phi \X & \longrightarrow & \X \\
\Big\downarrow && \mapdownright{\pi} \\
\X & \stackrel{\phi\circ\pi}{\longrightarrow} &  S
\end{array}\]
and
\item~$U\in\db{\X\times_\phi \X}$ is an object in the derived category of 
the product. 
\end{itemize}
There is a map~$\X\times_\phi\X\ra S$ with fibre~$X_s\times X_{\phi(s)}$ over 
$s\in S$. The (derived) restriction of~$U$ to this fibre gives a kernel
\[U_s=\pull y_s^*(U)\in \db{X_s\times X_{\phi(s)}}\]
where~$y_s\colon X_s\times X_{\phi(s)}\hookrightarrow\X\times_\phi \X$ is the 
inclusion. Hence a relative kernel defines a family of functors 
\[ \Psi_s=\Psi^{U_s}\colon \db{X_{\phi(s)}}\ra\db{X_s}.
\] 
In the present paper, a relative kernel~$(U,\phi)$ will be called
{\em invertible}, if for all~$s\in S$ the functor~$\Psi_s$ is a 
Fourier--Mukai functor. 
Let~$\daut{\X, S}$ be the group of invertible relative kernels up to isomorphism,
the {\em group of relative equivalences} of the family~$\X\ra S$. 
By construction, every element of the group~$\daut{\X, S}$ gives a family of 
Fourier--Mukai transforms over the base~$S$. 

The next statement is in some sense auxiliary, but it encompasses the point 
of view of the present article. Let~$X$ be a projective K3 surface or 
Calabi--Yau threefold. Let~$\pi\colon \X\ra S$ be a family of projective deformations 
of~$X$ over a polydisc~$S$, with~$\pi^{-1}(0)\cong X$ for~$0\in S$. Assume that 
the Kodaira--Spencer map
\[\psi\colon T_0S\ra H^1(X, \Theta_X)\] 
of the family is injective. 

Let~$U_0\in \db{X\times X}$ be an invertible kernel on~$X$ 
giving rise to a Fourier--Mukai functor~$\Psi=\Psi^{U_0}$ on~$X$. Using the Mukai 
map from the derived category to cohomology 
(see for example~\cite[Section 3.1]{ca}), 
there is an induced isomorphism
\[\psi\colon H^*(X, \C)\ra H^*(X, \C)\] 
preserving Hodge structures (in the sense of Mukai for the K3 case). 
In particular,~$H^{n,0}$ is preserved where~$n$
is the dimension of~$X$; so if~$\Omega\in H^0(X, \Omega_X^n)$ is a 
holomorphic top-form then its image~$\psi(\Omega)$ is also a holomorphic top-form
(a constant multiple of~$\Omega$). 

\begin{theorem}\label{defdirections} 
Assume that there is an invertible relative kernel~$(U,\phi)$ on~$\X\rightarrow S$
with~$\phi(0)=0$ extending~$U_0$. Then there is a commutative diagram
\[\begin{array}{ccc}
T_0(S) & \stackrel{d\phi|_0}\longrightarrow & T_0(S) \\
\mapdownright{\psi} && \mapdownright{\psi} \\
H^1(X, \Theta_X)&& H^1(X, \Theta_X)\\
\mapdownright{\wedge\Omega} && \mapdownright{\wedge\psi(\Omega)} \\
H^1(X,\Omega_X^{n-1})&&H^1(X,\Omega_X^{n-1}) \\
\Big\downarrow && \Big\downarrow\\
H^*(X, \C) & \stackrel{\psi}\longrightarrow & H^*(X, \C)
\end{array}\]
where the last vertical maps are the inclusions coming from Hodge theory. 
\label{cohaction}\end{theorem} 

This statement make look complicated, but it says something 
very simple. Suppose you have a Fourier--Mukai functor~$\Psi$ on~$X$. 
The action of~$\Psi$ on cohomology gives rise, via Hodge theory, to a 
map on the base of the local deformation space of~$X$. Then the only way to extend
$\Psi$ over a deformation family of~$X$ is to a relative kernel whose action~$\phi$
on the base is compatible with the map defined by~$\Psi$. In particular, 
unless~$\Psi$ acts trivially on the local deformation space, 
it will never extend to a family of autoequivalences
($\phi={\rm id}_S$) in a family of deformations of~$X$. 

\vspace{.05in}\noindent {\sc Proof of Theorem~\ref{cohaction}} \hspace{.05in}
Once the statement is properly formulated, the proof is not 
very difficult. Note that the family~$\Psi^{U_s}$ of Fourier--Mukai functors gives 
rise to an isomorphism of local systems 
$\oplus_nR^n(\phi\circ\pi)_*(\C_\X)\cong \oplus_nR^n\pi_*(\C_\X)$ on~$S$ (basically
just a continuous family of cohomology isomorphisms), which preserves Hodge 
filtrations. Now use the fact that the period map of the family is injective 
(since the Kodaira--Spencer map of~$\pi$ is, and~$X$ is Calabi--Yau), and unwind 
the definition of the derivative of the period map at~$0\in S$. 
\ethrm

\section{K3 surfaces with~$ADE$ configurations} 
\label{sec!K3}

Let~$\bar Y$ be a projective K3 surface with a du Val (rational double point) 
singularity at a point~$p\in\bar Y$ and no other singularities. 
Let~$g\colon Y\rightarrow \bar Y$ be its smooth K3 resolution with exceptional locus
$E= E_1\cup\ldots\cup E_r$. It is well known that each component~$E_i$ 
is a smooth rational curve of self-intersection~$-2$, hence 
it defines a reflection 
\begin{equation}
\label{defaction} r_i\colon \omega \mapsto\omega + \left(E_i\cdot\omega\right)E_i
\end{equation} 
on~$H^2(Y, \C)$. The intersection graph of the curves~$\{E_i\}$
is a Dynkin diagram~$\Delta$ of type~$ADE$, and as the notation suggests, the maps 
$r_j$ generate an action of the reflection group~$W_\Delta$ on~$H^2(Y, \C)$. 

\renewcommand{\labelenumi}{(\roman{enumi})}
\begin{proposition}\label{defs_of_Y} There exists a family~$e\colon \Y\rightarrow Z$ 
of projective deformations of~$e^{-1}(0)\cong Y$ over a complex polydisc 
$0\in Z$, with an action of the finite group~$W_\Delta$ on the base~$Z$, 
such that the following properties hold: 
\begin{enumerate} \item there is a proper subset~$Z_i\subset Z$ such that 
$s\in Z_i$ if and only if the fibre~$Y_s$ contains a smooth rational curve which 
is a deformation of~$E_i\in Y$; 
\item for every~$s\in Z$, 
there is a contraction morphism~$Y_s\rightarrow\bar Y_{i,s}$, which 
contracts the deformation of~$E_i$ in~$Y_s$ if~$s\in Z_i$ and is an 
isomorphism otherwise;
\item the fixed locus of~$r_i$ on~$Z$ equals~$Z_i$; and 
\item for~$w\in W_\Delta$ and~$s\in Z$, the fibres~$Y_s$,~$Y_{w(s)}$ are 
isomorphic.  
\end{enumerate}
\end{proposition} 
\proof This can be proved using the language of lattice-polarized 
K3 surfaces~\cite{dol}. Let~$M$ be the orthogonal complement of 
$\langle E_1,\ldots, E_n\rangle$ in the Picard group of~$Y$, or any sublattice 
thereof containing the cohomology class of an ample divisor on~$Y$; since~$\bar Y$ 
was assumed projective, such~$M$ exist. 
Consider the local moduli space~$\Y\ra Z$ of~$M$-polarized K3 surfaces~\cite{dol}
with central fibre~$Y=e^{-1}(0)$ for~$0\in Z$, a smooth family of projective K3 
surfaces. Since~$Z$ is small, the second cohomology~$H^2(Y_s, \Z)$ can be 
identified across the family.
The base~$Z$ is isomorphic, using the Kodaira--Spencer map, to a small
disc around the origin in~$N\otimes\C$, where~$N$ is the orthogonal complement of 
$M$ in~$\pic(Y)$. Since~$M$ does not include the class~$E_i$, 
$E_i\in H^2(Y_s, \Z)$ is algebraic (and represented by a rational curve) if
and only if~$s\in Z_i$ for a subvariety~$Z_i\subset Z$. It is easy to see that 
the~$W_\Delta$-action on~$H^2(Y,\C)$ preserves~$N\otimes\C$, and hence~$W_\Delta$
can be made act on~$Z$. The isomorphisms~$Y_s\cong Y_{w(s)}$ come from the 
Torelli theorem, since these surfaces have isomorphic Hodge structure. Finally 
the fact that~$Z_i$ is exactly the fixed locus of~$r_i$ is just chasing 
definitions. 
\ethrm

Next I~want to define relative kernels on~$\Y\rightarrow Z$, indexed by nodes
of the diagram~$\Delta$. By (ii) above, for a node~$i$ of~$\Delta$ and 
$s\in Z$ there is a
contraction~$Y_s\rightarrow \bar Y_{i,s}$ which contracts~$E_i$ if~$s\in Z_i$ 
and is an isomorphism otherwise. There is a diagram 
\[\begin{array}{ccccc}
&& \tilde Y_{i,s} \\
&\swarrow &&\searrow\\
Y_s&&&&Y_{r_i(s)}\\
&\searrow &&\swarrow\\
&&\bar Y_{i,s}
\end{array}\]
where~$\tilde Y_{i,s}$ is the fibre product of the two contractions. This fibre
product can be thought of as a subscheme of the product~$Y_s\times Y_{r_i(s)}$; 
it is the ``correspondence variety'' on the product (pairs of points  mapping
to the same image). 
If~$s\in Z\setminus Z_i$, then~$\tilde Y_{i,s}$ is simply the diagonal 
in~$Y_s\times Y_{r_i(s)}$ with respect to the isomorphism~$Y_s\cong Y_{r_i(s)}$. 
On the other hand, if~$s\in Z_i$, then~$E_i\subset Y_s$ is a rational curve, 
and~$\tilde Y_{i,s}$ has two components: 
one is the diagonal, and the other one is~$E_i\times E_i\cong \P^1\times \P_1$.
The components intersect along the diagonal~$\Delta_{E_i}$. 

In any case, set~$U_{i,s}=\O_{\tilde Y_{i,s}}\in\db{Y_s\times Y_{r_i(s)}}$
to be the (pushforward of the) structure sheaf of this correspondence subscheme. It 
is possible to show (see~\cite[Theorem 4.1]{sz_br} for the case of threefolds) that
the kernels~$U_{i,s}$ are restrictions to the fibres of a relative kernel 
$(U_i, r_i)$ on~$\Y\ra Z$. 

\begin{theorem} For every node~$i$ of~$\Delta$, the relative kernel 
$(U_i, r_i)$ is invertible: for~$s\in Z$, the kernel~$U_{i,s}$ defines a \fmt
\[ \Psi_{i,s}=\Psi^{U_{i,s}}\colon \db{Y_{r_i(s)}}\stackrel{\sim}\longrightarrow\db{Y_s}
\]
such that for a pair of nodes~$(i,j)$ of~$\Delta$, there is a isomorphism of functors
\begin{equation}\underbrace{\Psi_{i,s} \circ\Psi_{j,r_i(s)}\circ\ldots}_{m_{ij}} \cong \underbrace{\Psi_{j,s}\circ\Psi_{i,r_j(s)}\circ\ldots}_{m_{ij}}\colon\db{Y_{r_{ij}(s)}}\longrightarrow\db{Y_s}
\label{fmrel}\end{equation}
where
\[r_{ij}=\underbrace{r_i\circ r_j\circ\ldots}_{m_{ij}} \cong \underbrace{r_j\circ r_i\circ\ldots}_{m_{ij}}\in W_\Delta.
\]
Hence the derived category~$\db{Y}$ carries an action of the 
braid group~$B_\Delta$, and this action deforms to an action of~$B_\Delta$ by  
a family of derived equivalences over the deformation space~$\Y\ra Z$ of~$Y$.
\label{thm_K3} \end{theorem}
\proof The point~$s=0\in Z$ is fixed by all~$r_i$, and in this 
case the theorem is a re-statement of a special case of~\cite[Theorem 1.2]{st}.
In more detail, as proved in~\cite[Lemma 4.6]{sz_br}, for~$s=0\in Z$
the functors~$\Psi_{i,0}$ are just the twist functors of~\cite{st} 
with respect to the spherical sheaves~$\O_{E_i}(-1)$ on~$Y=Y_0$. 
The relations~(\ref{fmrel}) were proved in~\cite{st}. Hence mapping the 
braid group generator~$R_i$ to the autoequivalence~$\Psi_{i,0}$ defines 
an action of~$B_\Delta$ on~$\db{Y}$.

For arbitrary~$s\in Z$, the fact that~$\Psi_{i,s}$ is invertible is 
easy: if~$s\in Z_i$ then it is still a twist functor; otherwise it is the 
structure sheaf of the diagonal in~$Y_s\times Y_{r_i(s)}$ under the isomorphism 
$Y_s\cong Y_{r_i(s)}$, and hence clearly invertible. The relation~(\ref{fmrel}) 
can be proved using the method of~\cite{sz_br}, which does the more complicated 
case of threefolds. The point is that the kernels for the composites on both sides 
of the relation~(\ref{fmrel}) can be proved to be structure sheaves; 
for a general point~$s\in Z$ they are both isomorphic to the 
structure sheaf of the diagonal in~$Y_s\times Y_{r_{ij}(s)}$ under 
the isomorphism~$Y_s\cong Y_{r_{ij}(s)}$, and 
from this a specialization argument concludes that the two kernels are isomorphic 
everywhere. In particular, this gives an independent proof in this case of the 
braid relations on the central fibre~$Y$.  
\ethrm 

It is known from~\cite{aK3} that (for appropriate values of the K\"ahler form) 
Type II string theory on the surface~$Y$ exhibits enhanced gauge symmetry. 
The braid group action in Theorem~\ref{thm_K3} is a holomorphic 
shadow of this enhanced gauge symmetry: 
the derived category of~$Y$ has a braid group worth
of autoequivalences covering the Weyl group of the nonperturbative
gauge symmetry algebra, which deform to equivalences between different varieties
under a deformation of its complex structure. In other words,
at the enhanced gauge symmetry points the derived automorphism group of~$Y$
(the group of symmetries of the category of topological D-branes)
is larger than that of its deformations. 

I~next extend Theorem~\ref{thm_K3} and its interpretation as enhanced gauge symmetry
to gerby deformations, also known as nonzero~$B$-fields. 
I~take the most simple-minded definition, according to which 
the~$B$-field is a class~$B\in H^2(Y, \R/\Z)$. A~$B$-field can be used to twist the 
derived category of coherent sheaves of the K3 surface~$Y$ as follows. 
Consider the natural map 
\begin{equation}\label{bfieldmap} \delta\colon H^2(Y, \R/\Z)\rightarrow H^2(Y, \O_{X}^*)
\end{equation}
coming from the exponential sequence. The class 
$\beta=\delta(B)\in H^2(Y, \O_{X}^*)$ gives a ``gerbe'' on~$X$, and 
there is a notion of a sheaf over this gerbe (also called~$\beta$-twisted 
sheaf on~$Y$). One wants to define the ``derived category of~$\beta$-twisted sheaves 
on~$Y$'' with some finiteness condition. If the class~$B$ is torsion in 
$H^2(Y, \R/\Z)$, then the usual
notion of coherence generalizes, and one obtains~\cite{ca} a triangulated category 
$\db{Y, B}$ with properties very similar to those of~$\db{Y}$. 
In the general case there 
does not seem to be an accepted definition, though see~\cite[Remark 2.6]{ko} for 
discussion. The following statement is therefore formulated for the case
of torsion~$B$-fields; I~certainly expect it to hold in general. 

\begin{theorem} Let~$B\in H^2(Y, \Q/\Z)$ be a torsion~$B$-field. Then 
for every vertex~$i$ of~$\Delta$, there is a family of twisted \fmt s
\begin{equation}\label{twistfm} \Psi_{i,s,B}\colon \db{Y_{r_i(s)}, r_i(B)}\stackrel{\sim}\longrightarrow\db{Y_s, B}
\end{equation}
deforming the functor~$\Psi_{i,s,0}=\Psi_{i,s}$. Here
$W_\Delta$ acts on on~$H^2(Y, \R/\Z)$ via its action on~$H^2(Y, \R)$. 
\label{K3bfield}\end{theorem} 
\proof Let~$p_1, p_2$ denote the projections of~$Y_s\times Y_{r_i(s)}$
onto its factors. A twisted functor~(\ref{twistfm}) needs, 
by~\cite[Section 3.1]{ca}, a kernel 
\[V_{i,s}\in \db{Y_s\times Y_{r_i(s)}, p_2^*(r_i(B))-p_1^*(B)}\]
(note that I~am using additive notation for classes in cohomology with values in 
$\Q/\Z$). 

Recall the correspondence variety~$\tilde Y_{i,s}$ in~$Y_s\times Y_{r_i(s)}$ with 
respect to the~$i$-th contraction. The sheaf~$U_{i,s}$ was defined as the structure 
sheaf of this correspondence; more precisely, if 
$k\colon \tilde Y_{i,s}\hookrightarrow Y_s\times Y_{r_i(s)}$ is the inclusion, 
then~$U_{i,s}=k_*\O_{\tilde Y_{i,s}}$. 

Let \[\tilde B=p_2^*(r_i(B))-p_1^*(B).\] 
Note that by~\cite[Theorem 2.2.6]{ca}, there is a twisted pushforward functor
\[k_*\colon \db{\tilde Y_{i,s}, k^*(\tilde B)}\rightarrow \db{Y_s\times Y_{r_i(s)}, \tilde B}.\]
I~claim that the structure sheaf of the scheme~$\tilde Y_{i,s}$ is naturally a 
sheaf on~$\tilde Y_{i,s}$ over the gerbe defined by the class~$k^*(\tilde B)$. 
This implies that the kernel 
$U_{i,s}=k_*\O_{\tilde Y_{i,s}}$ can be 
thought of as a sheaf on~$Y_s\times Y_{r_i(s)}$
over the gerbe corresponding to~$\tilde B$
and hence it can be used to define the twisted functor~(\ref{twistfm}).

To prove the claim, I~distinguish two cases. First 
assume~$s\in Z_i$. Then~$E_i$ deforms to~$Y_s$ and
as I~said above,~$\tilde Y_{i,s}$ has two components: one is 
$E_i\times E_i$ and the other one is~$\Delta_{Y_s}$, the diagonal. It is enough to 
show that the structure sheaf of either component is a sheaf over the gerbe 
coming from~$k^*(\tilde B)$ restricted to that component. 
But one component~$E_i\times E_i$ is simply the
quadric surface, which has a trivial Brauer group and hence there is nothing to
prove. On the other component, 
$k^*(\tilde B)|_{\Delta_{Y_s}}=(B\cdot E_i)E_i$. Now the point
is that since~$s\in Z_i$,~$E_i$ is an algebraic class on~$Y_s$, hence the class
$k^*(\tilde B)|_{\Delta_{Y_s}}$ defines the trivial gerbe (see Remark~\ref{subtle} 
for the argument). Hence again, the structure sheaf is a sheaf over this gerbe! 

Next assume that~$s\in Z\setminus Z_i$. Then there is an isomorphism 
$Y_s\cong Y_{r_i(s)}$. It can be shown that this isomorphism induces the map 
$r_i$ on second cohomology. On the other hand,~$\tilde Y_{i,s}$ is in this 
case irreducible and isomorphic to the diagonal; moreover,~$\tilde B$ pulls back to 
the trivial gerbe over this diagonal. Hence the structure sheaf is again a sheaf
over the gerbe defined by~$k^*(\tilde B)$. 

The fact that the kernel~$U_{i,s}$ defines an equivalence of categories
can be proved using~\cite[Theorem 3.2.1]{ca}, which generalizes the criterion of 
Bridgeland~\cite[Theorems 5.1 and 5.4]{br}; I~omit the details. 
\ethrm

\begin{remark}\label{subtle}\rm The statement of Theorem~\ref{K3bfield}
involves a subtlety concerning the
$W_\Delta$-action on gerbes. On the central fibre~$Y$, all cohomology classes
$E_i$ are algebraic. On the other hand, the map~$H^2(Y, \R)\ra H^2(Y, \O_Y^*)$
factors through~$H^2(Y,\O_Y)$ and by Hodge theory, the image of 
$E_i\in H^2(Y, \R)$ in~$H^2(Y,\O_Y)$ is zero. 
This implies that~$B$ and~$r_i(B)$ give the same gerbe on~$Y$. However, for 
generic~$Y_s$ the classes~$E_i$ are transcendental, and~$B, r_i(B)$ are 
different gerbes. 
In fact, Theorem~\ref{K3bfield} should be complemented by a statement 
that there is no family of equivalences 
\[ \db{Y_{r_i(s)}, B}\rightarrow\db{Y_s, B}.
\]
The family of sheaves~$\{U_{i,s}\}$ is certainly not appropriate, 
since as the proof above shows, 
$\tilde B$ gives a nontrivial gerbe for~$s\in Z\setminus Z_i$ exactly because 
the class~$[E_i]$ is transcendental on~$Y_s$. 
Indeed I~expect that the only possible way to deform the equivalence
$\Psi_{i,s}$ in the~$B$-field direction is that compatible with its cohomology 
action; in other words, there is an analogue of Theorem~\ref{defdirections} 
for gerby deformations. I~have no idea how to prove this statement. 
\end{remark} 

I~wish to offer the following interpretation of Theorem~\ref{K3bfield}: 
Type II string theory 
on~$(Y, B)$ has enhanced gauge symmetry (for appropriate values of the
K\"ahler parameter) 
if and only if the derived category~$\db{Y, B}$ admits a set of twisted 
autoequivalences, which deform to twisted \fmt s between 
different points in moduli when the complex structure and~$B$-field parameters 
are deformed. Theorem~\ref{K3bfield},
together with Remark~\ref{subtle}, says that this is the case if and only if 
$r_i(B)=B$ for all~$i$, in other words if and only if~$E_i\cdot B=0$ for all
exceptional curves. Note that this condition on the~$B$-field 
is identical to that of~\cite[p.~4]{aK3}, found by an analysis involving
heterotic/Type II duality. 

\section{Calabi--Yau threefolds containing ruled surfaces} 
\label{sec!CY}

Let~$\bar X$ be a projective threefold with a curve of singularities
\[B={\rm Sing}(X)\hookrightarrow \bar X,\] 
such that along the 
curve~$\bar X$ has du Val singularities of uniform~$ADE$ type. 
The iterated blowup of the singular locus~$f\colon X\ra \bar X$ is a resolution of 
singularities. Locally over a point~$p\in B\in \bar X$, the fibre of~$f$ is a 
set of rational curves as before, intersecting according to the appropriate
$ADE$ type Dynkin diagram. However, globally there may be monodromy 
(see Figure~1):  
as~$p$ moves over the curve~$B$, the configuration of curves may be permuted 
according to a diagram symmetry of the Dynkin diagram. 

\begin{figure}[ht]
\centering
\input{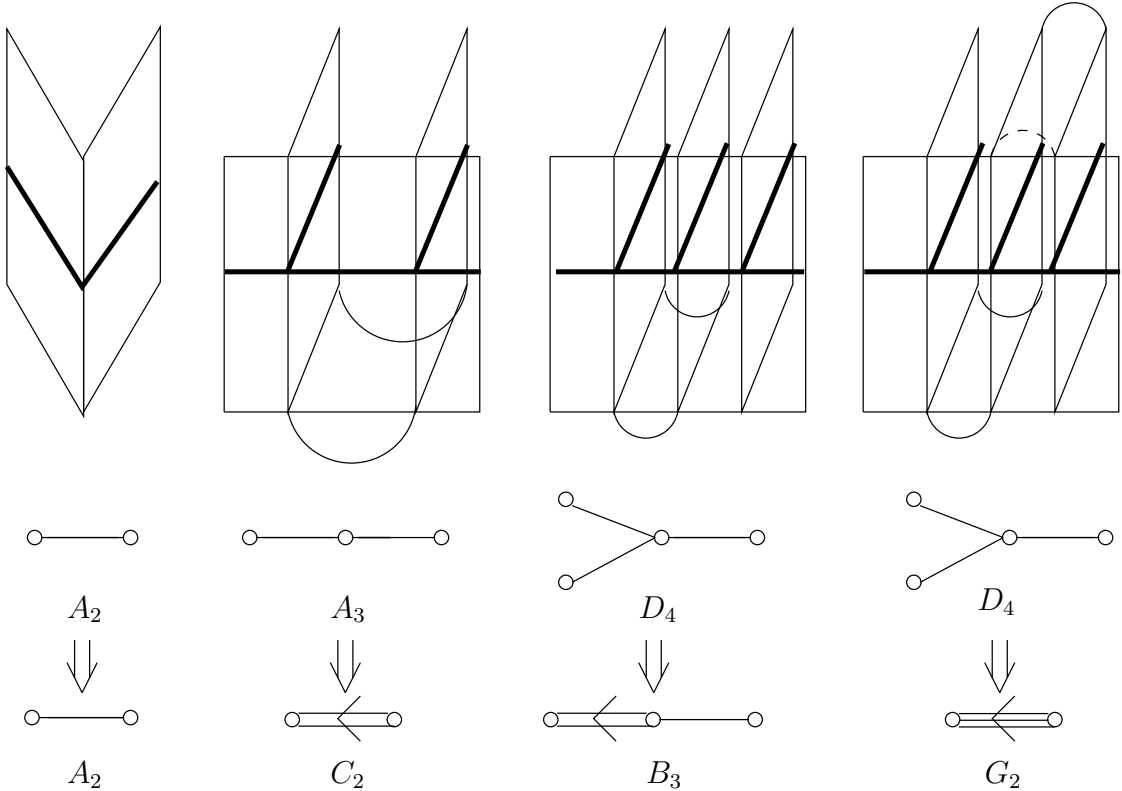}
\caption{Dynkin diagrams and configurations of surfaces} 
\end{figure}

It is well known that quotients of~$ADE$ Dynkin diagrams by (subgroups of) 
their automorphism groups are 
non-simply laced Dynkin diagrams in a well-defined sense. Concretely, the action 
of~$\Z/2$ on the diagrams~$A_{2n+1}$,~$D_n$ and~$E_6$ gives, respectively, 
the diagrams~$C_{n+1}$,~$B_{n-1}$ and~$F_4$, whereas the action of~$\Z/3$ and
the symmetric group on three letters leads to the diagram~$G_2$. The group  
$\Z/2$ also acts on the diagram~$A_{2n}$; this is a special case which I~exclude 
from consideration, though see Remark~\ref{specialcase} below. 

Globally therefore, the exceptional locus of~$f\colon X\ra \bar X$ 
consists of a set of smooth geometrically ruled surfaces 
$\{\pi_j\colon D_j\ra B_j\}$ intersecting in a Dynkin configuration~$\Delta$, 
which may or may not be simply laced. If~$\Delta$ is simply laced then each~$B_j$ is
isomorphic to~$B$, whereas in the general case each~$B_j$ is an umramified cover 
of~$B$ of the appropriate degree. 

As in the case of surfaces, I~want to describe some deformations of the
threefold~$X$. In the local case (when one restricts attention to a neighbourhood
of the exceptional surfaces), this problem 
is studied in detail in~\cite[Section 2]{sz_br}.
Globally there may be some obstructions to realizing all local deformations as
actual projective deformations of~$X$. In simple cases (see below) it can be 
checked that the deformations I~describe actually exist. The next proposition 
therefore should be considered a kind of ``ideal scenario'' statement. 

\begin{proposition} Let~$X$ be the Calabi--Yau threefold constructed above, 
with a set of exceptional surfaces~$\{\pi_j\colon D_j\ra B_j\}$ indexed by nodes of 
a Dynkin configuration~$\Delta$, which may or may not be simply laced. 
Assume that~$X$ has
good deformation theory. Then the universal family of (projective Calabi--Yau) 
deformations~$e\colon\X\ra S$ of~$X=e^{-1}(0)$ over a polydisc~$0\in S$ carries
an action of the reflection group~$W_\Delta$ on its base~$S$; moreover,
the following properties hold. 
\begin{enumerate} \item For every~$s\in S$, there is a contraction 
$f_s\colon X_s\ra \bar X_s$ deforming the contraction~$f$. 
\item There is an analytic subset~$S_j\subset S$ of codimension 
equal to the genus of~$B_j$, such that 
$s\in S_j$ if and only if the fibre~$X_s$ contains a smooth ruled surface 
in the exceptional locus of~$f_s$ which is a deformation of~$D_i\in X$.
\item The fixed locus of~$r_j$ on~$S$ equals~$S_j$.
\item For~$w\in W_\Delta$ and~$s\in S$, the fibres~$X_s$,~$X_{w(s)}$ are 
birational.
\end{enumerate}
Assume moreover that the genus~$g$ of~$B$ is at least one, 
and~$s\in S$ is a general point in the base. Then
\begin{enumerate}\addtocounter{enumi}{3}
\item the exceptional locus of~$X_s\ra\bar X_s$ consists of  
rational~$(-1,-1)$-curves, coming in sets of~$(2g-2)$
naturally indexed by positive roots of~$\Delta$.
\item For~$w\in W_\Delta$ and~$s\in S$, the birational map 
$X_s\dashrightarrow X_{w(s)}$ flops some of these curves.
\end{enumerate}
\label{defs_of_X}\end{proposition} 

Note that in the central fibre, the exceptional locus of~$f_s$ consists of a set 
of surfaces indexed by {\em simple} roots (nodes) of~$\Delta$. In the general fibre
(assuming genus at least two), the exceptional set of~$f_s$ is a set of curves
indexed by {\em positive} roots of~$\Delta$. 
Figure~2 illustrates the case~$\Delta=A_2, \,\,g=2$. 

\begin{figure}[ht]
\centering
\input{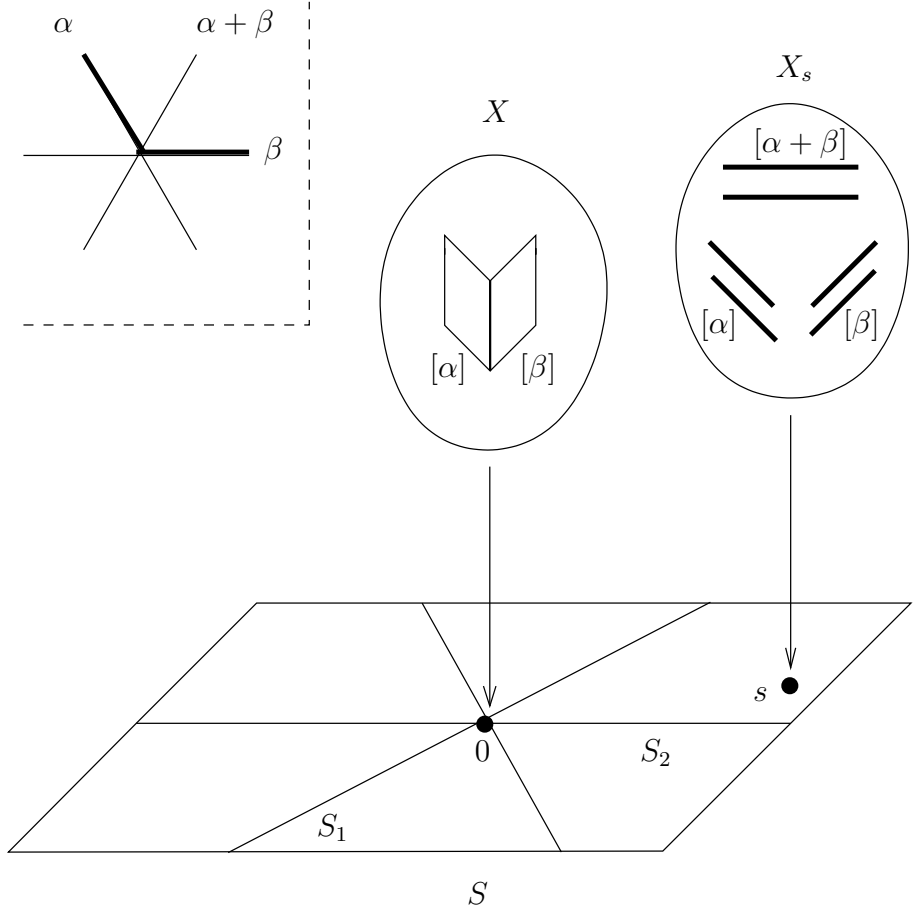}
\caption{The root system of~$A_2$ and exceptional loci for~$g=2$}
\end{figure}

Note also that the deformation theory of~$X$ is very different if the genus 
of~$B$ is zero. In that case, the~$W_\Delta$-action is trivial
($S_j=S$ for all~$j$ and hence every generator fixes~$S$) and the exceptional 
locus is always two-dimensional. For higher genus the~$W_\Delta$-action 
is non-trivial and the general exceptional locus is 
one-dimensional. The case~$g=1$ is also somewhat special: in that case, 
for general~$s\in S$, the contraction~$f_s\colon X_s\ra\bar X_s$ is an isomorphism, 
which is reminiscent of the surface case. 
This distinction is discussed further below.

The next statement is the exact analogue of Theorem~\ref{thm_K3}.

\begin{theorem}\label{thm_CY} For every node~$j$ of~$\Delta$, there is a family
of \fmt s
\[ \Psi_{j,s}\colon \db{X_{r_j(s)}}\stackrel{\sim}\longrightarrow\db{X_s}
\]
such that for a pair of nodes~$(i,j)$ of~$\Delta$, there is a isomorphism of 
functors
\begin{equation}\underbrace{\Psi_{i,s} \circ\Psi_{j,r_i(s)}\circ\ldots}_{m_{ij}} \cong \underbrace{\Psi_{j,s}\circ\Psi_{i,r_j(s)}\circ\ldots}_{m_{ij}}\colon\db{Y_{r_{ij}(s)}}\longrightarrow\db{Y_s}
\end{equation}
where
\[r_{ij}=\underbrace{r_i\circ r_j\circ\ldots}_{m_{ij}} \cong \underbrace{r_j\circ r_i\circ\ldots}_{m_{ij}}\in W_\Delta.
\]
Hence the derived category~$\db{X}$ carries an action of the 
braid group~$B_\Delta$, and this action deforms to an action of~$B_\Delta$ by  
a family of equivalences over the deformation space~$\X\ra S$ of~$X$.
\end{theorem}

\proof The proof, given in detail in~\cite[Section 4]{sz_br}, is similar to 
that of Theorem~\ref{thm_K3}. The individual functors~$U_{j,s}$ are defined 
using a diagram 
\[\begin{array}{ccccc}
X_s&&\dashrightarrow &&X_{r_j(s)}\\
&\searrow &&\swarrow\\
&&\bar X_s.
\end{array}\]
For $s\in S_i$, the functor turns out to be a special case of a functor 
written down by Horja in~\cite[(4.31)]{horja1}, and proved 
invertible in~\cite{horja2}. 
The proof of the braid relations relies, as before, on a specialization argument. 
\ethrm

According to \cite{aCY}, \cite{kmp}, \cite{akm} and 
references cited in these works, 
threefolds~$X$ of the above type (for suitable values of the K\"ahler form)
exhibit enhanced gauge symmetry. Theorem~\ref{thm_CY} is a holomorphic shadow of 
this symmetry: the derived category of~$X$ has a braid group worth of 
autoequivalences covering the Weyl group of the gauge algebra, which for genus at 
least one deforms away to a set of equivalences between different deformations.
In particular, the derived automorphism group of~$X$ is larger than generic 
at these enhanced symmetry points. 

It is interesting to consider the case when the curve~$B$ has genus zero. 
In this case, the projective threefold~$X$ has no deformations where the 
surfaces deform away. The braid group still acts on the derived category of
$X$, but it also acts as a set of derived autoequivalences on all deformations.
Hence nothing gets ``enhanced''. This phenomenon was also observed in the
physics literature: as explained in~\cite[p.2]{kmp}, enhanced gauge
symmetry needs that~$B$ is not rational; 
if~$B\cong\P^1$ then the symmetry is only present in the limit when the area 
of~$B$ goes to infinity~\cite{aCY}. The lack of deformations
is also an issue in the proof of the braid relations in~\cite{sz_br}; 
the proof proceeds via decomposing~$X$ locally into a union of two pieces
$X_1\cup X_2$, so that both contain ruled surfaces over the affine line
$\A^1$ and have enough deformations. Decomposing~$\P^1$ into a union
of two lines is here the mathematics equivalent to taking the area of 
the~$\P^1$ to infinity. 

\begin{examples} \rm 
Varieties~$\bar X$ with a curve of singularities of uniform 
type~$A_n$ can be found among hypersurfaces or complete intersections in weighted 
projectice spaces; compare for example~\cite{kmp}. The resolution~$X$ is then 
embedded in a (partial) resolution of the ambient space, typically with 
$n$ distinct divisors over the relevant singular locus; hence the configuration
in~$X$ is still of type~$A_n$. It can often be shown by concrete methods that the 
deformation theory of these threefolds is good in the sense needed for 
Proposition~\ref{defs_of_X} to hold. 
Such varieties can be systematically searched for and in low
codimension classified using the graded ring method pioneered by Reid; 
see the~$A_1$ case in~\cite{sz_wci} and the general case in~\cite{anita}.  

Just for amusement, I~proceed to give an example of a projective Calabi--Yau 
threefold~$X$ which contains a~$C_2$ configuration of surfaces, inspired 
by~\cite[Section 3]{akm}; to the best of my knowledge, this is the first
explicit example of this kind. Begin with an auxiliary variety
\[ \bar V = \left\{ \begin{array}{rcl} x_2^4&=&y_1y_2 \\ x_1^8 + x_2^8 + y_1^4 + y_2^4 + y_3^4 + z^2 &=&0\end{array}\right\} \subset \P^5[1,1,2,2,2,4]. 
\]
$\bar V$ is a degenerate degree~$(4,8)$ complete intersection Calabi--Yau threefold
in the indicated 
space. It can be checked by explicit computation that~$\bar V$ has three curves
of singularities, which are all elliptic. Along two of the curves at 
$\{x_1=x_2=y_1=0\}$ and~$\{x_1=x_2=y_2=0\}$,~$\bar V$ has generically~$A_1$
singulariries; this is a result of the identifications on the weighted projective
space. For a generic~$(4,8)$ complete intersection (which is simply an 
octic in $\P^4[1,1,2,2,2]$, since the degree four variable can be eliminated), 
there is one irreducible curve of~$A_1$ singularities, but in the special~$\bar V$  
this part of the singular locus
becomes reducible because of the first equation. The last curve is 
$\{x_2=y_1=y_2=0\}$, arising also because of the first equation; 
the singularity along the last curve is generically~$A_3$. 
The three curves all meet at the two points~$(0:0:0:0:1:\pm i)$ of the 
weighted projective space. A patient calculation shows that these points are 
also quotient singularities, under the group~$\Z/2\times \Z/4$ acting 
on~$\C^3$ by~$(-1,-1,1)\times (1,i,-i)$. 

Construct a particular crepant partial resolution~$V\ra \bar V$ in two steps. 
First perform the blowup of both intersection points according to the right hand 
arrow of the toric diagram Figure~3. This introduces
two exceptional divisors over the two points, and leaves behind three 
disjoint curves of singulatities of uniform type~$A_1$,~$A_1$ and~$A_3$ 
respectively, with no dissident points. Then blow up the two disjoint~$A_1$ curves
to get a Calabi--Yau threefold~$V$ with a single elliptic curve of uniform 
$A_3$ singularities. 

\begin{figure}[ht]
\centering
\input{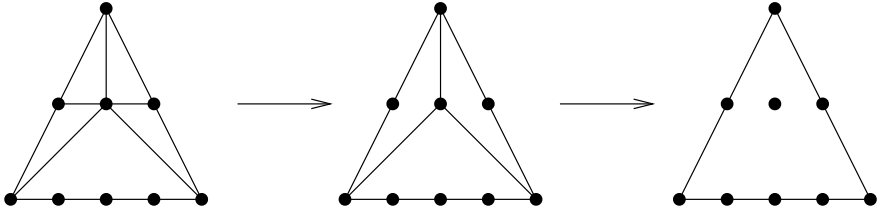}
\caption{The toric partial resolution of~$\C^3/(\Z/2\times \Z/4)$}
\end{figure}

Consider the action 
\[\iota\colon (x_1:x_2:y_1:y_2:y_3:z)\mapsto(x_1:(-x_2):y_2:y_1:(-y_3):(-z))\] 
on the weighted projective space. This action fixes~$\bar V$; since it 
interchanges the two~$A_1$ singular curves, it extends to the partial
resolution~$V$. Further,~$\iota$ acts by a free action on
the elliptic curve of~$A_3$ sigularities of~$V$; in the transverse coordinates
$x_2,y_1,y_2$ to this curve satisfying the relation~$x_2^4=y_1y_2$, the action 
interchanges~$y_1$ and~$y_2$, and maps~$x_2$ to~$-x_2$. A final 
check shows that~$\iota$ acts freely on~$\bar V$ and hence on~$V$. Thus letting
\[\bar X= V/\langle\iota\rangle,\]
the projective Calabi--Yau threefold~$\bar X$ has an elliptic curve of~$A_3$ 
singularities and is smooth otherwise; moreover, the local coordinates along this
curve undergo~$\Z/2$ monodromy. Hence its Calabi--Yau resolution~$X\ra \bar X$ 
contains a~$C_2$ configuration of exceptional surfaces ruled over 
elliptic curves. 
\end{examples} 

\begin{remark} \rm The braid group action on the derived category gives rise 
to actions on even and odd cohomology, using the Mukai map. 
The action on odd cohomology~$H^3(X,\C)$ leads, as discussed 
in Proposition~\ref{defdirections}, to a Weyl group action 
on the tangent space to the deformation space, in a
compatible fashion with the way the derived equivalences deform. 
There is also an induced Weyl group action on the Picard group. 
Some of these actions were described before; eg.~\cite{kmp} has a 
symmetric group action in the case of Type A, both on the Picard group and
the deformation space. The action of the braid group on the 
derived category explains all these actions in a uniform way. 
\end{remark}

\begin{remark}\label{specialcase}
\rm The case of monodromy~$\Z/2$ acting on the Dynkin diagram 
$A_{2n}$ has been excluded from consideration all along. This case has caused 
considerable headache also in the physics literature~\cite[Section 4]{akm}. 
In this case, the exceptional divisors~$D_i$ of~$f\colon X\ra\bar X$
are still indexed by vertices of a kind of quotient quiver,
the~$A_n$-quiver with a marked vertex at one end corresponding 
to the adjacent~$\Z/2$-orbit of vertices of~$A_{2n}$.  
However, the marked node corresponds to a singular 
exceptional surface. It is an irreducible non-normal surface 
$\pi_n\colon D_n\ra B$ whose double locus is a 
section and whose fibre over any point $b\in B$ is a line pair. 
I~do not know whether there exists an autoequivalence $\Psi_n$ 
of~$\db{X}$ corresponding to this surface, but I~suspect that the answer is 
yes; this is a contracting~$EZ$-configuration in the sense of Horja~\cite{horja2}, 
with singular~$E$. 

\begin{figure}[ht]
\centering
\input{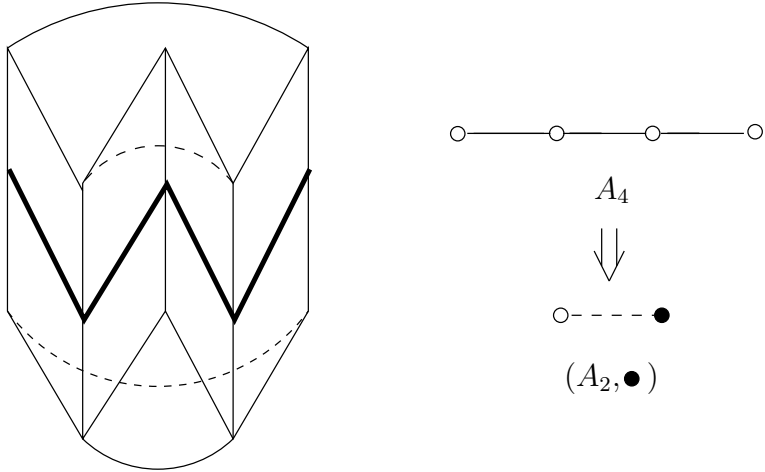}
\caption{The $A_4$ configuration with $\Z/2$ monodromy}
\end{figure}

The Main Assertion of \cite[Section 4]{akm}, supported by various arguments 
including the analysis of the matter spectrum, claims that the enhanced 
gauge symmetry is~${\mathfrak s \mathfrak p}(n)$, or in the language of the present 
paper, of type~$C_n$. The point of view exposed in 
this paper gives additional support to this claim. Namely, 
the derived category of~$\db{X}$ is acted on by the autoequivalences 
$\Psi_1, \ldots, \Psi_{n-1}$ coming from the smooth ruled surfaces, as well as
the hypothetical autoequivalence~$\Psi_n$; the question is what are the relations. 
One can make an educated guess based on the following argument. 

In the singular threefold~$B\subset\bar X$, take a small quasiprojective surface 
$\bar Y\subset\bar X$ intersecting~$B$ once transversally. 
Let~$Y\ra \bar Y$ be its resolution 
in~$X$, with exceptional curves~$E_1, \ldots, E_{2n}\subset Y$. 
Set~$\E_i={\O_{E_i}(-1)}\in\db{Y}$ for~$i=1,\ldots, 2n$.
The functors~$\Psi_i$ can be restricted to Fourier--Mukai 
functors on~$Y$ (compare~\cite[Proof of Theorem 4.5]{sz_br}).  
The functor~$\Psi_i$ for~$1\leq i\leq n-1$ restricts in fact to the composite 
of two twist functors~$T_{\E_i}$ and~$T_{\E_{2n+1-i}}$. On the other hand,
by~\cite{st}, the twist functors~$\{T_{\E_i} \colon1\leq i\leq 2n\}$ generate 
the braid group~$B_{A_{2n}}$ acting on the derived category of $Y$. Moreover, the 
monodromy~$\Z/2$ acts on this braid group, 
mapping~$T_{\E_i}\mapsto T_{\E_{2n+1-i}}$ for~$i=1,\ldots ,n$. 
The guess I~want to make is that
the functors~$\Psi_1, \ldots, \Psi_n$ satisfy the relations of the fixed 
subgroup~$(B_{A_{2n}})^{\Z/2}$. By a result in algebra~\cite{michel}, this fixed 
subgroup is generated by the composites~$T_{\E_i}\circ T_{\E_{2n+1-i}}$
(note these commute) for~$i=1,\ldots , n-1$ and a final element 
$T_{\E_n}\circ T_{E_{n+1}}\circ T_{\E_n}$ (note these braid), and  
the group generated by these elements is the braid group corresponding to the Weyl
group~$(W_{A_{2n}})^{\Z/2}$. This latter group can be checked by a direct argument
to be isomorphic to the Weyl group of the diagram~$C_n$. 

Hence the conjectural answer is that~$X$ has a set of derived equivalences
$\Psi_1, \ldots, \Psi_n$ satisfying the braid relations of the 
Dynkin diagram~$C_n$. In other words,~$X$ has enhanced gauge symmetry of type~$C_n$
(or~${\mathfrak s \mathfrak p}(n)$). 
\end{remark}

\begin{remark}\rm To conclude this section, I~remark that as opposed to 
the case of dimension two, the braid group actions of~\cite{st} can never be
interpreted as enhanced gauge symmetry in dimension three. The reason is the
following: it can easily be checked that if~$\E$ is a spherical object in the 
sense of \cite{st}, then the corresponding twist functor acts on cohomology by
$\alpha\mapsto\alpha + \langle\ch(\E), \alpha\rangle \ch(\E)$, where 
$\langle,\rangle$ is a linear combination of intersection forms on cohomology. 
However,~$\ch(\E)$ only has even components, hence the action of the twist 
functor on odd cohomology and so on~$H^1(X,\Theta_X)$ is trivial. In particular, 
by Theorem~\ref{defdirections}, a twist functor always deforms to all deformations 
as an autoequivalence in dimension three, and hence it can never be part of an 
``enhanced'' action. 
\end{remark} 

\section{Elliptic fibrations and braid groups of affine type} 
\label{sec!aff}

Let~$\sigma\colon X\ra S$ be an elliptic fibration of a projective threefold~$X$. 
Assume that there is a smooth component~$C\subset S$ of the discriminant locus of 
$\sigma$, over which the fibres of~$\sigma$ are of uniform Kodaira type~$I_n$
with~$n>2$, $I^\star_n$, $II^*$, $III^*$ or $IV^*$. These are the fibre types 
corresponding to the affine diagrams $\tilde A_{n-1}\,\, (n>2)$, $\tilde D_{n+4}$, 
$\tilde E_6$, $\tilde E_7$ and $\tilde E_8$. In~$X$, the rational curves in 
the fibres over~$p\in C$ undergo monodromy, and trace out ruled 
surfaces~$\pi_j\colon D_j\rightarrow C_j$. Assume that in the type $\tilde A_{n-1}$ 
case the monodromy is not transitive, and in the type $\tilde D_4$ case it does not 
act transitively on the outer vertices. Then the global intersections of the
exceptional surfaces are described by an affine Dynkin diagram $\tilde\Delta$, 
which is the original $\tilde A\tilde D\tilde E$ diagram for trivial 
monodromy and a quotient non-simply laced $\tilde B\tilde C\tilde G\tilde F$ type
diagram otherwise.
The diagram~$\tilde\Delta$ gives rise to a braid group~$B_{\tilde\Delta}$, 
with one generator for every node of~$\tilde\Delta$ and one 
(braid) relation for every pair of nodes as dictated by the labels of the
diagram~$\tilde\Delta$.

\begin{figure}[ht]
\centering
\input{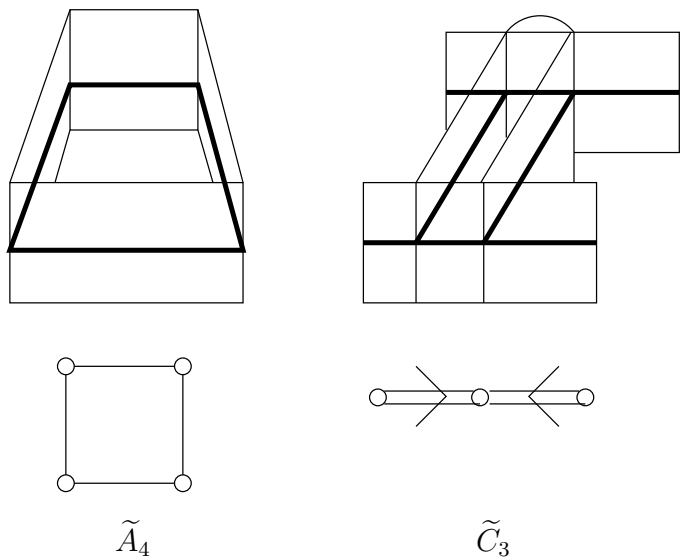}
\caption{Some ruled surface configurations in elliptic fibrations} 
\end{figure}

\begin{theorem} The braid group~$B_{\tilde\Delta}$ of affine 
type acts on the derived category~$\db{X}$. 
\end{theorem} 
\proof The ruled surfaces~$D_j\ra C_j$ give rise to \fmt s~$\Psi_j$ on~$X$
as before. The proof of a single braid relation only concerns two surfaces and
the functors defined by them. Under the assumptions made, every pair of surfaces
forms an~$A_1\times A_1$,~$B_2$ or~$G_2$ configuration. Moreover, the computation
of the composed functors can be restricted to a small neighbourhood of these two 
surfaces. Hence the proof of~\cite{sz_br} applies. 
\ethrm

Enhanced gauge symmetry for threefolds with elliptic fibrations has been 
discussed in the context of F-theory compactifications; 
see~\cite{mv}, \cite{akm}, \cite{mg} and references therein.

\section{Mirror symplectomorphisms?}
\label{sec!mir}

The paper~\cite{st}, a direct predecessor of the present work, is directly 
motivated by mirror symmetry. Namely, the original motivation of that 
paper was to find the mirrors of certain symplectomorphisms of symplectic 
manifolds~$(M^{2n},\omega)$, {\em Dehn twists} in Lagrangian spheres 
$S^n\subset M$. The twist functors in spherical objects are natural candidates
for the mirrors of Dehn twists. 

As discussed in~\cite{horja1}, \cite{sz} and \cite{ahk}, 
the derived equivalences studied in this paper, arising from ruled surfaces 
collapsing to curves in~$X$, are mirror to certain diffeomorphisms of the 
mirror manifold, arising as monodromy transformations around certain 
boundary components of the complex moduli space of the mirror~$M$. 
These diffeomorphisms are symplectomorphisms of~$(M^{2n},\omega)$
for special values of the symplectic form~$\omega$. 
It would be of interest to find a direct symplectic geometric
construction of these diffeomorphisms. It is tempting to speculate that they 
are given by some kind of twisting with respect to a fibered submanifold of~$M$, 
just as the \fmt s of~$X$ are constructed from the ruled surfaces. 
\cite{kklm} begins the topological study of the mirrors of some explicit
Calabi--Yau manifolds~$X$ containing a single ruled surface;  
the situation appears to be quite intricate. It would also 
be interesting to see whether in appropriate cases the braid 
relations~(\ref{braidrelations}) can be proved for these symplectomorphisms.

\vspace{0.2in} 

\noindent {\small \sc Department of Mathematics, Universiteit Utrecht, PO. Box 80010, NL-3508 TA Utrecht, The Netherlands}

\noindent and

\noindent {\small \sc Alfr\'ed R\'enyi Institute of Mathematics, Hungarian Academy of Sciences, PO. Box 127, H-1364 Budapest, Hungary}

\noindent {\small E-mail address: \tt szendroi@math.uu.nl}
\end{document}